\documentclass{amsart}

\usepackage[russian,english]{babel}

\theoremstyle{plain}

\numberwithin{equation}{section}

\title[On $\pi$-extensions of the semigroup $\mathbb{Z}_+$]{On $\pi$-extensions of the semigroup $\mathbb{Z}_+$}
\author{T.A. Grigoryan, \ E.V. Lipacheva, \ V.H. Tepoyan}

\begin{document}

\maketitle

\begin{center}
\textit{ Kazan State Power Engineering University }
\end{center}

\textbf{Abstract. }{ We study inverse $\pi$-extensions of the semigroup $\mathbb{Z}_+$. It is shown that $\pi$-extension of the semigroup 
$\mathbb{Z}_+$ is inverse, iff its $\pi$-extension coincides with $\pi(\mathbb{Z}_+)$. The existence of a non-inverse 
$\pi$-extension for semigroup $\mathbb{Z}_+$ is proved. }



\textbf{1. Introduction.} In his well known work [1] Coburn proved that all the isometric representations of the semigroup
 $\mathbb{Z}_+$ of non-negative integers generate canonically isomorphic $C^*$-algebras. This theorem was generalized by many authors
 to a larger class of semigroups. Douglas [2] showed the same for the semigroup of positive cone of real numbers
 $\mathbb{R}$. Murphy proved this theorem for the positive cones of abelian groups with order. On the other hand, Murphy [3] and 
Jang [4] have shown that this theorem is not true for the semigroup $\mathbb{Z}_+\backslash \{1\}$. The isometric representations
 with commuting range projections of the semigroup $\mathbb{Z}_+\backslash \{1\}$ has been studied by Raeburn and Vittadello [5]. 

We introduce the notion of $\pi$-extension of the semigroup of non-negative integers $\mathbb{Z}_+$ (see definition 2.1), and study
 the properties of $\pi$-extensions of the semigroup $\mathbb{Z}_+$. Also the concept of the inverse $\pi$-extension of the semigroup
 $\mathbb{Z}_+$ is introduced in definition 5.1. We prove, that if $\pi$ is an irreducible representation, then there is no
 non-trivial inverse $\pi$-extension for this semigroup. In case $\pi$ is reducible, there exists a non-inverse $\pi$-extensions. 
On the other hand we show that for any isometric representation of $\mathbb{Z}_+$ there always exists $\pi$-extension. 

The authors are sincerely grateful to S. A. Grigoryan for useful discussions and for valuable advice.

\textbf{2. Preliminaries.} Consider an isometric representation of the semigroup $\mathbb{Z}_+$:
$$\pi:\mathbb{Z}_+\rightarrow B(H),$$
where $B(H)$ is a set of all bounded linear operators on Hilbert space $H$. Denote by $Is(H)$ the semigroup of all isometric operators
 in $B(H)$.

\textit{Definition 2.1.} We call $M\subset Is(H)$ a $\pi$-extension of the semigroup $\mathbb{Z}_+$, if:\\
	1. $\pi(\mathbb{Z}_+)\subset M$;\\
	2. $\pi(i)T=T\pi(i)$ for $T$ in $M$ and $i$ in $\mathbb{Z}_+$.

The following irreducible representation is considered throughout this paper, unless the opposite is mentioned: 
$$\pi:\mathbb{Z}_+\rightarrow B(H^2),$$
where $H^2(S^{1},d\mu)$ is the Hardy space of square-integrable complex-valued functions on the unit circle $S^{1}$ by Haar measure
 $\mu$, and with spectrum in $\mathbb{Z}_+$.

The operator $\pi(n)$ is the multiplicative operator of multiplication by the function $e^{in\theta}$, i.e.:
$$\pi(n):H^2\rightarrow H^2 \mbox{ and } \pi(n)f(e^{i\theta})=e^{in\theta}f(e^{i\theta}).$$

The orthonormal system of functions $1$, $e^{i\theta}$, $e^{i2\theta}$, $\dots$ form a basis in $H^2$, and the operator $\pi(1)$ is
 the shift operator on this basis:
$$\pi(1)e^{in\theta}=e^{i(n+1)\theta}.$$
Therefore, the $C^*$-subalgebra of the algebra $B(H^2)$ generated by the operators $\pi(1)$ and $\pi^*(1)$ is a Toeplitz algebra.

\textbf{3. Inverse Representations.} \emph{An inverse semigroup} $P$ is a semigroup, such that each element $x$ has a unique 
\emph{inverse} element $x^*$ satisfying
$$xx^*x=x, \ x^*xx^*=x^*.$$

We denote by $\Delta_{S}$ the set of all isometric representations of the semigroup $S$. For ${\pi}\in\Delta_{S}$ define $S^{\pi}$
 to be the semigroup generated by operators ${\pi}(i)$ and ${\pi}^{*}(i)$, where $i\in{S}$.

\textit{Definition 3.1.} We call the representation ${\pi}\in{\Delta_{S}}$ inverse, if $S^{\pi}$ is an inverse semigroup.

\emph{The regular isometric representation} is a map ${\pi}:S\rightarrow{B(l^2(S))}$, $i\mapsto {\pi}(i)$, defined by following relation:
\[
({\pi}(i)f)(j) =
\begin{cases}
f(k), & \text{if $j=i+k$ for some $k\in S$;} \\
0, & \text{otherwise.}
\end{cases}
\]
\textbf{\textit{T h e o r e m 3.1.}} The regular isometric representation of the semigroup $S$ is inverse (see [6, 7]).

Now we give an example of non-inverse representation.
 
Let $\pi:\mathbb{Z}_+\rightarrow B(H^2)$ be the representation of the semigroup $\mathbb{Z}_+$ described in the section 2, i.e.
 $\pi(n)$ is the multiplicative operator of multiplication by the function $e^{in\theta}$.

Every inner function $\Phi(z)$ defines an isometric multiplicative operator $T_{\Phi}$:
$$T_{\Phi}f=\Phi f.$$

\textbf{\textit{T h e o r e m 3.2.}} Let $\widetilde{\pi}:\mathbb{Z}_+\times\mathbb{Z}_+\rightarrow B(H^2)$ be a representation, which maps
 $(n,0)\mapsto e^{in\theta}$ and $(0,m)\mapsto \Phi^m$, where $\Phi$ is an arbitrary inner function, not in semigroup 
$\{e^{in\theta}\}_{n=0}^{\infty}$. Then $\widetilde{\pi}$ is a non-inverse representation, i.e.  $(\mathbb{Z}_+\times\mathbb{Z}_+)^{\widetilde{\pi}}$ is a non-inverse semigroup.

\textbf{4. $\pi$-extension of the semigroup $\mathbb{Z}_+$.}

\textbf{\textit{L e m m a 4.1.}} Every isometric operator in $\pi$-extension of the semigroup $\mathbb{Z}_+$ can be
 represented through a single inner function.

Let us denote by $C^{*}_{\pi}(\mathbb{Z}_+)$ the $C^{*}$-algebra generated by the isometric representation $\pi$, described in Section 
2. Let also $C^*_{\pi}(M)$ be the $C^{*}$-algebra, generated by a semigroup $M\subset Is(H^2)$.

If $M$ is a $\pi$-extension of the semigroup $\mathbb{Z}_+$, then by Lemma 4.1 for each isometric operator $T\in M$ there exists 
 a unique inner function $\Phi$, such that the operator $T$ is a multiplication operator by $\Phi$. Define 
$$M^{'}=\{\Phi; \ T_{\Phi}\in M\}.$$

\newpage

\textbf{\textit{T h e o r e m 4.1.}} Let $M$ be the $\pi$-extension of the semigroup $\mathbb{Z}_+$. Then the following
 conditions are equivalent:

1. $C^{*}_{\pi}(\mathbb{Z}_+)=C_{\pi}^*(M);$
	
2. $M^{'}$ is a subsemigroup of the semigroup of finite Blaschke products.

\textbf{5. Inverse $\pi$-extension.} We denote by $\mathbb{Z}_+^{\pi}$ the involutive semigroup generated by $\pi(\mathbb{Z}_+)$ and $\pi(\mathbb{Z}_+)^*$.
 Let $M$ be the $\pi$-extension of the semigroup $\mathbb{Z}_+$. Denote by $\mathcal{M}^*$ the semigroup generated by $M$ and $M^*$. 
\textit{Definition 5.1.} We call the $\pi$-extension of the semigroup $\mathbb{Z}_+$ inverse, if $\mathcal{M}^*$  is an inverse semigroup.

Let $\pi:\mathbb{Z}_+\rightarrow B(H^2)$ be the representation of the semigroup $\mathbb{Z}_+$, described in Section 2. Then the
 following result is true. 

\textbf{\textit{T h e o r e m 5.1.}} $\mathcal{M}^*$ is inverse iff $\mathcal{M}^*=\mathbb{Z}_+^{\pi}$.

Consider an arbitrary isometric representation $\pi:\mathbb{Z}_+\rightarrow B(H)$. Denote $H_0=\ker{\pi^*(1)}$. It is clear that $H_0$
 is a Hilbert subspace of $H$.

\textbf{\textit{T h e o r e m 5.2.}} Let $\pi:\mathbb{Z}_+\rightarrow B(H)$ be an isometric representation of the semigroup
 $\mathbb{Z}_+$ such that the subspace $H_0=\ker{\pi^*(1)}$ is not one dimensional. Then there exists an inverse $\pi$-extension $M$
 of the semigroup $\mathbb{Z}_+$ such that $\mathbb{Z}_+^{\pi}$ is a proper involutive subsemigroup of the involutive semigroup $\mathcal{M}^*$.

\begin{center}
R E F E R E N C E S
\end{center}

1. \textbf{Coburn L.A.} \textit{The C*-algebra Generated by an Isometry}. // Bull. Amer. Math. Soc., 1967, v. 73, p. 722-726.

2. \textbf{Douglas R.G.} \textit{On the C*-algebra of a One-parameter Semigroup of Isometries}. // Acta Math., 1972, v. 128, p. 143-152.

3. \textbf{Murphy G.J.} \textit{Ordered Groups and Toeplitz algebras}. // J. Operator Theory, 1987, v. 18, p. 303-326.

4. \textbf{Jang S.Y.} \textit{Generalized Toeplitz Algebras of a Certain Non-amenable Semigroup}. // Bull. Korean Math. Soc., 2006, v. 43, p. 333-341.

5. \textbf{Raeburn I., Vittadello S.T.} \textit{The Isometric Representation Theory of a Perforated Semigroup}. // J. Operator Theory, 2009, v. 62, p. 357-370.

6. \textbf{Aukhadiev M.A., Tepoyan V.H.} \textit{Isometric Representations of Totally Ordered Semigroups}. // Lobachevskii Journal of Mathematics, 2012, v. 33, p. 239-243.

7. \textbf{Grigoryan S.A., Salakhutdinov A.F.} \textit{C*-algebras Generated by Cancellative Semigroups}. // Sibirsk. Mat. Zh., 2010, v. 51, p. 16-25.

{\small \vspace{\baselineskip}\hrule \vspace{3pt}
\par
{\bf Grigoryan Tamara Anatolevna} -- {Kazan State Power Engineering University}
\par
E-mail: {\it tkhorkova@gmail.com} }
{\small \vspace{\baselineskip}\hrule \vspace{3pt}
\par
{\bf Lipacheva Ekaterina Vladimirovna} -- Kazan State Power Engineering University
\par
E-mail: {\it elipacheva@gmail.com} }
{\small \vspace{\baselineskip}\hrule \vspace{3pt}
\par
{\bf Tepoyan Vardan Hakobi} -- Kazan State Power Engineering University
\par
E-mail: {\it tepoyan.math@gmail.com} }

\end{document}